\documentclass[oneside,english,reqno]{amsart}
\usepackage{amssymb}
\usepackage{graphics}

\newcommand\R{\ensuremath{\mathbb{R}}}

 \theoremstyle{plain}
 \newtheorem*{lem*}{Lemma}
 \theoremstyle{plain}
 \newtheorem{prop}{Proposition}
 \theoremstyle{remark}
 \newtheorem*{rem*}{Remark}
 \newtheorem{example}{Example}
 \newtheorem*{prob*}{Problem}
 \theoremstyle{plain}
 
 \newtheorem{lemma}{Lemma}
 \def\trace{{\rm trace}}

\def\E{{\sf E}}

\def\Pr{{\sf P}}
\def\T{^{\sf T}}
\def\eps{\varepsilon}

\def\th@nopoint{
  \thm@headpunct{} 
  \itshape 
} \theoremstyle{nopoint}

\begin{document}

\title{Waiting for a bat to fly by (in polynomial time)}

\author{Itai Benjamini}

\author{Gady Kozma}

\author{L\'aszl\'o Lov\'asz}

\author{Dan Romik}

\author{G\'abor Tardos}

\begin{abstract}
We observe returns of a simple random walk on a finite graph to a
fixed node, and would like to infer properties of the graph, in
particular properties of the spectrum of the transition matrix. This
is not possible in general, but at least the eigenvalues can be
recovered under fairly general conditions, e.g. when the graph has a
node-transitive automorphism group. The main result is that by
observing polynomially many returns, it is possible to estimate the
spectral gap of such a graph up to a constant factor.
\end{abstract}
\maketitle

\section{Introduction}

A spelunker has an accident in the cave. His lamp goes out, he cannot
move, all he can hear is a bat flying by every now and then on its
random flight around the cave. What can he learn about the shape of
the cave?

In other words: What can we learn about the structure of a finite
graph using only information obtained by observing the returns of a
random walk on the graph to this node?

Let $G=(V,E)$ be a connected simple graph with $n=|V|>1$ vertices, and let
$r\in V$ be a fixed node. Let $w_0=r, w_1, w_2, \dots$ be the steps
of a simple random walk on $G$ starting from $r$. Assume that we
observe the {\it return time sequence}, the infinite sequence of
(random) times $0<T_1<T_2< \dots$ when the walk visits $r$.
Alternatively this can be described as a sequence $a_1, a_2, a_3, ...
$ of bits, where $a_i=1$ if the walk is at $r$ at time $i$, $0$
otherwise. Note that $T_2-T_1, T_3-T_2,\dots$ are independent samples
from the same distribution as $T_1$, which we call the {\it return
distribution} of $G$ to $r$.

We say that a parameter $p(G,r)$ of the graph $G$ and root $r$ can be
reconstructed (from the return time sequence), if for every two
rooted graphs $(G,r)$ and $(G',r')$ for which the return time
sequence has the same distribution, we have $p(G,r)=p(G',r')$.

Which graph parameters can be reconstructed from the return time
sequence? There is a trivial way to construct different graphs with
the same return sequence: take two isomorphic copies and glue them
together at the root. Sometimes it makes sense to assume that we also
know the degree $d(r)$ of the root. In this case, we can reconstruct
the number of edges through
\begin{equation}\label{EDGENUM}
|E|=d(r)\E(T_1)/2.
\end{equation}
If the graph is regular, then we can reconstruct the number of nodes:
\begin{equation}\label{NODENUM}
n=|V|=\E(T_1).
\end{equation}

Another trivial example is to observe if all the numbers $T_i$ are
even. This is so if the graph is bipartite, and it happens with
probability 0 otherwise.

A natural candidate for a reconstructible quantity is the spectrum of
the transition matrix $M$ of the random walk on $G$. Let
$\lambda_1=1, \lambda_2, ..., \lambda_n$ be the eigenvalues of $M$,
arranged in decreasing order. Bipartiteness is equivalent to saying
that $\lambda_n=-1$.


We are going to show by a simple example that the spectrum is not
reconstructible in general. On the other hand, we show that if
$\lambda$ is an eigenvalue of $G$ which has an eigenvector $v\in\R^V$
such that $v_r\not=0$, then $\lambda$ is reconstructible.
We note that the {\it multiplicity} of $\lambda$ is not necessarily
reconstructible.

A special case where the eigenvector condition above is satisfied for
all eigenvalues is when $G$ is node-transitive. We don't know whether
in this case the multiplicities are reconstructible.

Of particular interest is the issue of {\it efficient
reconstruction}, by which we mean observing a polynomial (or expected
polynomial) number of returns. We consider this question in the case
of the \emph{spectral gap} $\tau= 1- \lambda_2$. Assuming the graph
is node transitive, we describe a procedure to estimate $\tau$ up to
a constant factor, using just polynomially many (in $n$) of the first
values of the $T_i$. We give an example of a graph where the spectral
gap cannot be recovered \emph{at all} from observations made at one
particular node.


This question was first mentioned, together with other related
problems, in \cite{BL}. Another related work is that of Feige
\cite{F} which presents a randomized space-efficient algorithm that
determines whether a graph is connected. His method uses return times
of random walks to estimate the size of connected components.

\section{Examples}

\begin{example}\label{TREES}
Consider the two trees in Figure \ref{2-TREES}. The distribution of
the return time to the root is the same in both trees (see later).
The eigenvalues of the tree on the left are
\[
1, \sqrt{3}/2,  \sqrt{6}/4, 0, 0, 0, 0, 0, -\sqrt{6}/4, -\sqrt{3}/2,
-1,
\]
while the eigenvalues of the tree on the right are
\[
1, \sqrt{3}/2, \sqrt{3}/2, \sqrt{6}/4, 0, 0, 0, -\sqrt{6}/4,
-\sqrt{3}/2, -\sqrt{3}/2, -1.
\]
Note that the eigenvalues are the same, but their multiplicities are
different.

\begin{figure}[htb]
\includegraphics*{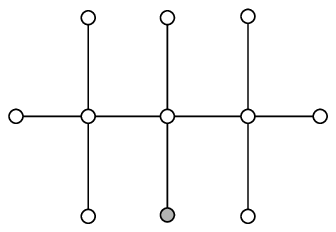}\includegraphics*{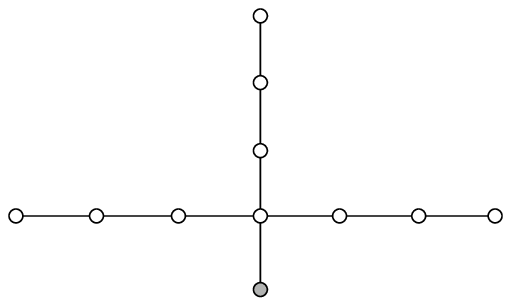}
\caption{Two trees with the same return times but different spectra}
\label{2-TREES}
\end{figure}
\end{example}

\bigskip

\begin{example}\label{EXPAND}
Let $T$ be a tree in which all internal nodes have degree $d+1$ and
which has a ``root'' $r$ such that all leaves are at distance $h$
from the root. We construct a graph $G$ by adding a $d$-regular graph
on the leaves.

For a fixed $h$ and $d$, all graphs obtained this way are
$(d+1)$-regular graphs, and the distribution of the return time to
the root is the same in all such graphs. On the other hand, graphs
obtained this way can have very different properties. If we add an
expander on the leaves, the graph $G$ will be an expander.  (Recall
that G is a $c$-expander iff $|\partial S| > c |S|$ for every non
empty set of vertices $S$ with $|S| < |G|/2$. For background on
expanders and spectral gap see e.g.\ \cite{Lu}.) If we connect
``twin'' leaves to each other, and also match up ``cousins'' to get
$d$ new edges at each node, then for $h>2$ the root will be a
cutpoint. For expanders, the eigenvalue gap $\lambda_1-\lambda_2$ is
bounded from below by a positive function of $d$, while for the
graphs with cutpoints in the middle the eigenvalue gap tends to 0 as
$h\to\infty$.
\end{example}

\section{Preparation: some algebra and generating functions}

\subsection{Return probabilities and eigenvalues}

Denote by $P_k(x,y)$ the probability that a simple random walk on $G$
starting at $x\in V$ will be at $y\in V$ at time $k$. Clearly
\begin{equation} \label{HIT}
P_k(x,y) = e_x\T  M^k e_y.
\end{equation}
Here $M$ is not symmetric, but we can consider the symmetrized
matrix $N=DMD^{-1}$, where $D$ is a diagonal matrix with the positive
numbers $\sqrt{d(i)}$ in the diagonal. The matrix $N$ has the
same eigenvalues as $M$, and so we have
\begin{equation}\label{RET-EXP}
P_k(r,r) = \sum_{i=1}^n f_i(r)^2 \lambda_i^k,
\end{equation}
where $f_1, f_2, ... , f_n$ is an orthonormal basis of eigenfunctions
of $N$ corresponding to the eigenvalues $\lambda_1,
\lambda_2,...,\lambda_n$.

We note that if the graph is node-transitive, then the value
$P_k(r,r)$ is the same for all $r$, and hence by averaging
(\ref{RET-EXP}) we get the simpler formula
\begin{equation}\label{RET-TRACE}
P_k(r,r) = \frac{1}{n}\trace(M^k) = \frac{1}{n} \sum_{i=1}^n
\lambda_i^k.
\end{equation}

At some point, it will be convenient to consider the {\it lazy
version} of our chain, i.e., the Markov chain with transition matrix
$M'=(1/2)(I+M)$ (before doing a step, we flip a coin to decide if we
want to move at all). The observer can easily pretend that he or she
is watching the lazy version of the chain: after each step, he flips
a coin in quick succession until he tosses a head, and advances his
watch by the number of coinflips. The distribution after $k$ lazy
steps is easy to compute from (\ref{HIT}):
\begin{equation}\label{HIT-LAZY}
P'_k(x,y)=2^{-k} e_x\T (I+M)^k e_y = 2^{-k}\sum_{j=0}^k \binom{k}{
j} e_x\T M^j e_y = 2^{-k}\sum_{j=0}^k \binom{k}{ j} P_j(x,y).
\end{equation}

The main advantage of the lazy chain is that its eigenvalues are
nonnegative. Furthermore, for a lazy chain we have
\[
\lambda_2+\dots+\lambda_n= \trace(M)-1 = \frac{n}{ 2}-1,
\]
and hence $\lambda_2\ge 1/3$ if $n\ge 4$.

%

\subsection{The generating function of return times}

Let us introduce the generating function
\begin{equation}\label{RET-GEN}
f(t)=\sum_{k=0}^{\infty} P_k(r,r)t^k = \sum_{i=1}^n f_i(r)^2
\frac{1}{1-t\lambda_i}.
\end{equation}
There are several other useful expressions for $f(t)$; for example,
we get from (\ref{HIT}) that
\[
f(t)= e_r\T (I-tM)^{-1} e_r,
\]
and expressing this in terms of determinants, we get
\begin{equation}\label{RET-DET}
f(t)= \frac{\det(I'-tM')}{\det(I-tM)},
\end{equation}
where $M'$ is the matrix obtained from $M$ by deleting the row and
column corresponding to the root, and $I'$ is the $(n-1)\times(n-1)$
identity matrix.

It will be convenient to do a little algebraic manipulation. The
reciprocal of this function is also an interesting generating
function:
\begin{equation}\label{FIRSTRET-GEN}
\frac{1}{f(t)}=1-\sum_{k=1}^{\infty} s_kt^k,
\end{equation}
where $s_k=\Pr(T_1=k)$ is the probability that the first return to
the root occurs at the $k$-th step. This function has a root at
$t=1$, so it makes sense to divide by $1-t$, to get the analytic
function
\begin{equation}\label{SURV-GEN}
\frac{1}{(1-t)f(t)}=\sum_{k=0}^{\infty} z_kt^k,
\end{equation}
where
\[
z_k=1-\sum_{j\le k} s_k=\sum_{j>k} s_k
\]
is the probability that the random walk does not return to the root
during the first $k$ steps.


\section{Reconstructing nondegenerate eigenvalues}

It is these formulas which form the basis of learning about the
spectrum of $G$ from the visiting times of the random walk at $x$,
since $P_k(r,r)$ is determined by the distribution of return times,
and can be easily estimated from the visiting times (see section
\ref{EFFICIENT}). We call an eigenvalue of $M$ {\it nondegenerate} if at least
one of the corresponding eigenfunctions $f(x)$ satisfies $f(r)\not= 0$. One
can see from (\ref{RET-EXP}) that the non zero nondegenerate eigenvalues are
determined by the distribution of return times. Using $\sum_{i=1}^nf_i(r)^2=1$
for the orthonormal basis $f_i$ we conclude that whether zero is a
nondegenerate eigenvalue of $M$ is also determined. The return time
distribution determines $f(t)$ and this can also be used to find the
nondegenerate eigenvalues: the poles of $f(t)$ are exactly the reciprocals
of the non zero, nondegenerate eigenvalues of
$M$. Zero is a nondegenerate eigenvalue if and only if
$\lim_{t\to\infty}f(t)>0$. Then we get

\begin{prop}\label{NONZERO}
If two rooted graphs have the same return time distribution, then
they have the same nondegenerate eigenvalues.
\end{prop}

Let us remark that if $G$ has a node-transitive automorphism group,
then {\it every eigenvalue of $M$ is nondegenerate.} Indeed, every
eigenvalue has an eigenvector, which does not vanish at some node;
by node-transitivity, it also has an eigenvector that does not vanish
at the root.

Let us also remark that the {\it multiplicity} of a nondegenerate
eigenvalue is not uniquely determined: $0$ is a nondegenerate
eigenvalue of both trees in Example \ref{TREES}, but it has different
multiplicities in the two. Furthermore, degenerate eigenvalues are
not determined by the return times: the second largest eigenvalues
of the transition matrices of the two $(d+1)$-regular graphs
constructed in Example \ref{EXPAND} are different. It follows from
Proposition \ref{NONZERO} that at least for the second graph, the
second largest eigenvalue is degenerate.


\section{Trees}\label{TREESEC}


We want to put Example \ref{TREES} in broader context. For trees, we
can simplify the generating function a bit: Since trees are
bipartite, we have $z_{2k}=z_{2k+1}$, and hence it makes sense to
divide by $t+1$ and then substitute $x=t^2$. It will be convenient to
scale by the degree of the root, and to work with the function
\begin{equation}\label{DEF-h}
h_G(x) =d(r)\sum_{k=0}^{\infty} z_{2k} x^k=
\frac{d(r)}{(1-x)f(\sqrt{x})}.
\end{equation}
It is easy to see that we did not lose any information here: we have
$h_{G_1}(x)=h_{G_2}(x)$ for two
trees $G_1$ and $G_2$ if and
only if they have the same return time distribution and their roots have the
same degree.

For a rooted tree with a single edge, $h_G(x)=1$. If a rooted tree
$G$ is obtained by gluing together the roots of two rooted trees
$G_1$ and $G_2$, then
\begin{equation}\label{GLUE-TREE}
h_G(x)=h_{G_1}(x)+h_{G_2}(x).
\end{equation}
This is easily seen by conditioning on which tree the random walk
starts in. Furthermore, if we attach a new leaf $r'$ to the root $r$
of a tree $G$ and make this the root to get a new rooted tree $G'$,
then
\begin{equation}\label{ADD-ROOT}
h_{G'}(x)=\frac{1+h_G(x)}{1+(1-x)h_G(x)}.
\end{equation}
To see this, consider a walk on $G'$ starting at $r'$, and the
probability $z'_{2k}$ that it does not return to $r'$ in the first
$2k$ steps ($k\ge 1$). The first step leads to $r$; the second step
has to use a different edge, which has a probability of
$d(r)/(d(r)+1)$. We can view the walk now as a random walk on $G$
until it returns to $r$. The probability that this happens after $2j$
steps is $z_{2j-2}-z_{2j}$. If $j\ge k$ then the walk will certainly
not return to $r'$ in the first $2k$ steps. If $j<k$, then we can
think of the situation as just having made a step from $r'$, and so
the probability that we don't return to $r'$ in the next $2k-2j-1$
steps is $z'_{2k-2j}$. Hence we get the equation
\[
z'_{2k}= \frac{d(r)}{d(r)+1}\left( z_{2k-2} +
\sum_{j=1}^{k-1}(z_{2j-2}-z_{2j})z_{2k-2j}\right.
\]
Multiplying by $x^k$ and summing over all $k\ge 0$, we get
(\ref{ADD-ROOT}).

These formulas can be verified from the definition of $z_k$. They
imply that $h_G$ is a rational function with integral coefficients.
They also provide us with a fast way to compute $h_G$, and through
this, to verify that the two trees in Example \ref{TREES} have the
same return distribution. But we can get more, a way to generate many
such pairs.

Suppose that we find a linear dependence between functions $h_G$ for
various trees $G$. This can be written as
\[
a_1h_{G_1}+\cdots+a_kh_{G_k} = b_1 h_{G_1'} +\cdots+b_m h_{G_m'}
\]
with some positive integers $a_1,\dots,a_k,b_1,\dots,b_m$. Now if we
glue together the roots of $a_1$ copies of $G_1$, $\dots$, $a_k$
copies of $G_k$ to get $G$, and the roots of $b_1$ copies of $G_1'$,
$\dots$, $b_m$ copies of $G_m'$ to get $G'$, then by
(\ref{GLUE-TREE}) we'll have
\[
h_G(x)= h_{G'}(x).
\]
We can add a new root to both if we prefer to have an example rooted
at a leaf.

Obviously, we only need to look for trees rooted at leaves. To find
such linear dependencies, it is natural to find trees for which
$h_G(x)$ is ``simple'', namely the ratio of two linear functions, and
then find three with a common denominator.  A general example is a
tree $G=G_{a,b}$ of height $3$, where the neighbor of the root has
degree $a$ and has $a-1$ neighbors of degree $b$. We can allow the
degenerate cases $b=1$ (when $G$ is a star rooted at a leaf) and
$a=1$ (when $G$ is a single edge). It is easy to compute that
\footnote{Are these the only trees for which $h_G$ has rational
numerator and denominator? Can one say anything about quadratic? What
about depth 4?}
\[
h_G=\frac{ab-(b-1)x}{ab-(ab-1)x}.
\]


So if we fix a $k$ which is not a prime, and consider trees
$G=G_{a,b}$ with $ab=k$, they all have the same denominator
$k-(k-1)x$, and so for any three of them their functions $h_G$ will
be linearly dependent. The simplest choice is $k=4$, when we get the
trees $G_{1,4}$ (a single edge), $G_{2,2}$ (a path of length 3) and
$G_{4,1}$ (a 4-star). Simple computation shows that
\[
h_{G_{1,4}}-3h_{G_{2,2}}+2h_{G_{4,1}}=0.
\]
Gluing these together as described above, and adding a new root for
good measure, gives the two trees in Example \ref{TREES}.

Using (\ref{RET-DET}) and (\ref{DEF-h}), it is not hard to see that the roots
of the numerator of $h_G(x)$ are the squared reciprocals of the
nondegenerate non zero eigenvalues of $G$, except for the trivial
nondegenerate eigenvalues $\pm 1$. The multiplicities, as we have
seen, are not necessarily determined by $h_G$.


\medskip
\begin{rem*} In the special trees constructed above, the squareroots
of the root of the denominator are exactly the degenerate eigenvalues
of $G$. We don't know if this is always so. An interesting open
question seems to be whether the degenerate eigenvalues are
reconstructible for trees.
\end{rem*}

\section{Effective reconstruction}\label{EFFICIENT}

In the previous section, we assumed that the exact distribution of
the return time is known, which is the same as saying that we can
observe the random walk forever. In this section we are concerned
with determining quantities after observing a polynomial number of
returns.

\subsection{Estimating return probabilities}


We show that we can estimate $P_k(r,r)$ from the observation of
polynomially many return times. Fix $k$ and observe the returns
$T_1,T_2,\ldots$ until the first $T_{i_1}$ with $T_{i_1}\ge k$; call this
period an {\it experiment}. Call the experiment {\it successful} if
$T_{i_1}=k$. The probability that an experiment is successful is
$P_k(r,r)$. Note that observing the next $k$ steps and then until the first
return (i.e., $T_{i_1+1},\ldots, T_{i_2}$ with the smallest $i_2$ such that
$T_{i_2}\ge T_{i_1}+k$) is an independent experiment.

So we have a sequence of independent events with the same
probability $p=P_k(r,r)$, and we want to estimate $p$. By standard
results, observing $p\eps^{-2}\delta^{-1}$ of them, the relative
frequency will be closer than $\eps$ to $p$ with probability
$1-\delta$.

The amount of time a particular trial takes is a random variable,
whose expectation is $k$ plus the time it takes to get back to $r$
after $k$ steps. This can be bounded by the maximum hitting time
between nodes, which is $O(n^3)$. Summing up,

\begin{prop}\label{PKFIND}
In an expected time of $O((k+n^3)\eps^2\delta^{-1})$ we can compute
an estimate of $P_k(r,r)$ which is within an (additive) error of
$\eps$ with probability $1-\delta$.
\end{prop}


\subsection{Reconstructing the eigenvalue gap}

We restrict our attention to node-transitive graphs, in which case we
can use the trace formula (\ref{RET-TRACE}). We can use
(\ref{NODENUM}) to reconstruct the number of nodes $n$. Furthermore,
we assume that the chain is lazy, so that its eigenvalues are
nonnegative, and their sum is $n/2$.

For a lazy chain, $P_k(r,r)$ tends to $1/n$ monotone decreasing.
Furthermore, (\ref{RET-TRACE}) implies that setting
\[
q_k=P_k(r,r)-\frac{1}{n},
\]
we have
\[
n q_{k+1}=\sum_{i=2}^n \lambda_i^{k+1} \ge
\frac{1}{n-1}\left(\sum_{i=2}^n \lambda_i\right)\left(\sum_{i=2}^n
\lambda_i^k\right) = \frac{1}{n-1} (\trace(M)-1) n q_k,
\]
and hence
\begin{equation}\label{QDECREASE}
q_{k+1} \ge \frac{1}{3} q_k
\end{equation}
for $n\ge4$ (which we assume without loss of generality).

We can try to compute recursively $\lambda_1 = 1$ and
$$
\lambda_i = \lim_{k\to\infty} \left[P_k(r,r)-\sum_{j=1}^{i-1}
\frac{\lambda_j^k} n \right]^{1/k}.
$$
This, however, does not seem to give an effective means of estimating
$\lambda_i$ in polynomial time. But to estimate at least the
eigenvalue gap $\tau=1-\lambda_2$ we can use the following fact.


\begin{lemma}\label{GAP}
We have
\begin{equation}\label{TAUBOUND}
\left(1+\frac{\ln n}{\ln q_k}\right)(1-q_k^{1/k}) \le \tau \le
1-q_k^{1/k}.
\end{equation}
\end{lemma}

It is not hard to see that these bounds imply the weaker but more
informative bounds
\begin{equation}\label{LOGBOUND}
\frac{\ln (1/q_k)}{k-\ln (1/q_k)} \le \tau \le \frac{\ln (n/q_k)}{k}.
\end{equation}

\begin{proof}
From (\ref{RET-TRACE}),
$$
P_k(r,r)=\frac{1}{n}+\sum_{i=2}^n \frac{\lambda_i^k}{n},
$$
and hence
$$
\frac{\lambda_2^k}{n} \le \sum_{i=2}^n \frac{\lambda_i^k}{n} = q_k
\le \lambda_2^k.
$$
Thus
\[
1-(nq_k)^{1/k} \le \tau \le 1-q_k^{1/k}.
\]
Using the elementary inequality
\[
\frac{1-x}{1-y}\le \frac{\ln x}{\ln y}
\]
valid for $0<x<y<1$, (\ref{TAUBOUND}) follows.
\end{proof}

Let $c>1$. It follows that if we find an integer $k>0$ such that
$q_k<1/n^c$, then $1-q_k^{1/k}$ is an estimate for the eigenvalue gap
$\tau$ which is within a factor of $1+1/c$ to the true value. But of
course we don't know $q_k$ exactly, only with an additive error: by
proposition \ref{PKFIND}, we can estimate $q_k$ in polynomial time with an
additive error less than (say) $\eps/n^c$, with high probability. So
to get valuable information, we need to find a value of $k$ for
which $q_k
>\eps/n^c$.

It is well known that the eigenvalue gap of a graph with $n$ nodes is
at least $1/n^2$, so we get that for $k\ge K_0=(c+1)n^2\ln n$,
\[
q_k \le n \left(1-\frac{1}{ n^2}\right)^k < ne^{-k/n^2}<
\frac{1}{n^c}.
\]

Applying Proposition \ref{PKFIND}, we can compute an approximation $Q_k$ of
$q_k$ that is within an additive error of $\eps/(8n^c)$ with
probability $\delta/(\log_2 K_0)$. By binary search, we can find a $k$ in the
interval $[0,K_0]$ for which $Q_k \le 1/n^c$ but $Q_{k-1}>1/n^c$.


\begin{prop}
For the value of $k$ computed above, $1-Q_k^{1/k}$ is within a
factor of $1\pm\eps$ of $\tau$ with probability at least $1-\delta$.
\end{prop}

\begin{proof}
With large probability, we have
\[
|q_m-Q_m|<\frac{\eps}{8n^c}
\]
for all $m$ for which we compute $Q_m$, in particular for $m=k-1$
and $m=k$. Using (\ref{QDECREASE}),
\[
q_k \ge \frac{1}{3} q_{k-1}  \ge
\frac{1}{3}\left(Q_{k-1}-\frac{\eps}{8n^c}\right) \ge \frac{1}{4n^c},
\]
and also
\begin{equation}\label{QHATQ0}
Q_k \ge q_k - \frac{\eps}{8n^c} \ge (1-\frac{\eps}{2})q_k.
\end{equation}
Similarly,
\begin{equation}\label{QHATQ1}
Q_k\le (1+\frac{\eps}{2})q_k.
\end{equation}
We claim that
\begin{equation}\label{QHATQ2}
1-\frac{\eps}{2} \le \frac{1-Q_k^{1/k}}{1- q_k^{1/k}}\le
1+\frac{\eps}{2}.
\end{equation}
To show the upper bound, we may assume that $Q_k\le q_k$. Then using
(\ref{QHATQ1}),
\[
\frac{1-Q_k^{1/k}}{1- q_k^{1/k}} \le \frac{\ln Q_k}{\ln q_k} \le
\frac{\ln((1-\frac{\eps}{2})q_k)}{\ln
q_k}=1+\frac{\ln(1-\frac{\eps}{2})}{\ln
q_k}<1-\ln(1-\frac{\eps}{2})\le 1+\frac{\eps}{2}.
\]
The lower bound in (\ref{QHATQ2}) follows similarly. Hence by Lemma
\ref{GAP},
\[
\tau \ge 1-q_k^{1/k} \ge (1-\eps)(1-Q_k^{1/k}),
\]
and
\begin{align*}
\tau &\le 1-\left(\frac{q_k}{n}\right)^{1/k} \le \left(1+\frac{\ln n}
{\ln (1/q_k)} \right) (1-q_k^{1/k})\\
&\le\left(1+\frac{1}{c}\right)\left(1+\frac{\eps}{2}\right)(1-Q_k^{1/k})
\le (1+\eps)(1-Q_k^{1/k}).
\end{align*}
\end{proof}

\section{Concluding remarks}

1. We can estimate for every node-transitive graph, by similar means,
the value $1-\max(\lambda_2,|\lambda_n|)$, which governs the mixing
time of the chain. The trick is to consider the matrix $M^2$ instead
of $M$, i.e., observe the chain only every other step. A little care
is in order, since this new chain may not be connected; but by
node-transitivity, its eigenvalue gap is the eigenvalue gap of the
component containing the observation node.

\medskip

2. The second moment of the first return time also has some more
direct meaning. Let $H(\pi,r)$ denote the expected number of steps
before a random walk starting from the stationary distribution hits
the root $r$. Then it is not hard to show using that the walk is
close to stationary at a far away time that
\[
H(\pi,r)=\frac{\E(T_1^2)}{2\E(T_1)}-\frac12.
\]
It is not clear whether any of the higher moments have any direct
combinatorial significance.

\medskip

3. Here are a couple of related problems.

\medskip

\noindent {\bf Problem:} Let $G$ be a connected graph of size $n$. We
label the vertices randomly by $m(n)$ colors and observed the colors
as they are visited by a simple random walk random walk: after each
step, the walker tells you ``now I'm at red'', ``now at blue'', and
so on. How many colors are needed in order to recover the shape of G
a.s. from this sequence of colors?

\medskip

\noindent {\bf Problem:} Consider an $n$-node connected graph. Take
$n$ particles labeled $1,...,n$. In a configuration, there is one
particle at each node. The interchange process introduced in
\cite{al} is the following continuous time Markov chain on
configurations: For each edge $(i,j)$ at rate $1$ the particles at
$i$ and $j$ interchanged. Assume you observed the restriction of the
interchange process to a fixed node, what graph properties can be
recovered? Obviously you get more information than in the case
discussed in the paper, which corresponds to noticing only one of the
particles. But is it really possible to use this information to
discover more about the graph?

\end{document}